\documentclass[11pt]{article}
\usepackage{amsmath, amsthm}
\usepackage{amssymb}
\usepackage{graphicx}
\newtheorem{theorem}{Theorem}

\newtheorem{definition}[theorem]{Definition}

\newtheorem{notation}[theorem]{Notation}
\newtheorem{remark}[theorem]{Remark}

\begin{document}

\title{On the weights of simple paths in weighted complete graphs}
\author{  Elena Rubei }
\date{\hspace*{1cm}}
\maketitle

\vspace{-1.1cm}

{\small 
{\bf Address:}
Dipartimento di Matematica ``U. Dini'', 
viale Morgagni 67/A,
50134  Firenze, Italia 

{\bf
E-mail address:}
 rubei@math.unifi.it}

\bigskip

\smallskip

\def\thefootnote{}
\footnotetext{ \hspace*{-0.36cm}
{\bf 2000 Mathematical Subject Classification: 05C22} 

{\bf Key words:} graphs,  weights of graphs}

{\small {\bf Abstract.} 
Consider a weighted graph $G$ with $n$ vertices, numbered by the set 
$\{1,...,n\}$.  For any path $p$ in $G$, we call $w_G(p) $ 
 the sum of the weights of the edges of the path and 
we define the multiset 
$${\cal D}_{i,j} (G) = \{w_G(p) | \; p \; simple \; 
path \; between \; i\, and \;j \}$$
We establish a criterion to say when, given a multisubset of $ \mathbb{R}$,
there exists a weighted complete graph $G$ 
such that the multisubset is equal to ${\cal D}_{i,j} (G)$ for some $i,j$ vertices of 
$G$. 

Besides we establish a criterion to say when, given for
 any $i, j \in \{1,...,n\}$ a multisubset of $ \mathbb{R}$,
${\cal D}_{i,j} $, there exists a weighted complete 
graph $G$ with vertices     $ \{1,...,n\}$  such that 
${\cal D}_{i,j} (G)= {\cal D}_{i,j} $ for any $i,j$.}

\bigskip

%matrice incompleta =coefficient matrix 
%matrice completa = augmented matrix

\section{Introduction}

The problem of the realization of metrics or, more generally,
positive symmetric matrices by graphs has a very rich literature.
The problem can be described shortly as follows.
 Given a weighted graph $G$  
(that is a graph such that every 
edge is endowed with a  real number, which we  call 
 the weight of the edge), for any path $p$ in $G$,
 we call $w_G(p) $ 
 the sum of the weights of the edges of the path and, 
for any $ i,j $ vertices of $G$, we define 
$$D_{i,j} (G) = min \{w_G(p) | \; p\; path \; between \; i\; and \; j\}$$ 
Given positive real numbers $D_{i,j}$ for any $i,j$ in a finite set $X$, we 
can wonder whether 
  there exists a positive-weighted graph $G$ whose set of vertices 
contains $X$ and such that $D_{i,j} (G) = D_{i,j}$ for any $ i, j \in X$. 

In  \cite{HY} Hakimi and Yau proved that a positive symmetric matrix  
$(D_{i,j})_{i,j}$
is realizable by a weighted graph if and only if it is a metric.
 
In \cite{B} Buneman established a criterion to see if a metric on a finite set
$X$ can be 
realized by a  positive-weighted tree with $X$ as set of leaves
(a partial result in this direction
 had already been obtained in \cite{SimP}). Precisely he proved 
that given a metric  on $\{1,...,n\}$,
$(D_{i,j})_{i,j \in \{ 1,...,n \}}$, there exists  
 a positive-weighted tree $T$  with  $\{1,...,n\}$ as set of leaves 
such that $ D_{i,j}= D_{i,j}(T)$  
if and only if, for all $i,j,k,h  \in \{1,...,n\}$,
the maximum of $\{D_{i,j} + D_{k,h},D_{i,k} + D_{j,h},D_{i,h} + D_{k,j}
 \}$ is attained at least twice.

In  \cite{B-S} Bandelt and Steel proved a result, analogous to
 Buneman's one, for not necessarily positive weighted trees:
for any set of real numbers $\{D_{i,j}\}_{i,j \in \{1,...,n\}}$,
 there exists a weighted tree $T$
such that $ D_{i,j} (T)= D_{i,j}$  for any $i,j  \in  \{1,...,n\}$  if 
and only  if, for any $a,b,c,d \in  \{1,...,n\}$,  we have that at least two
 among 
 $ D_{a,b} + D_{c,d},\;\;D_{a,c} + D_{b,d},\;\; D_{a,d} + D_{b,c}$
are equal.

The problem of reconstructing weighted 
trees from data involving the distances  
between the leaves has several applications, such as phylogenetics and 
some algorithms to reconstruct trees from the data $\{D_{i,j}\}$ have been 
proposed (among them we quote neighbour-joining method, invented by
 Saitou and Nei in 1987, see \cite{NS}, \cite{SK}).

Obviously the problems of realization of symmetric matrices by graphs and of 
reconstructing the weighted graphs from the ``distances'' 
between the vertices may have some applications, 
also in the case  the weights are not all positive or all negative.
Imagine that a particle, by going through an edge 
of a graph, gets or looses some substance (as much as the weight of the
 edge). 
If we know how much the substance of this particle varies by going from a 
vertex $i$ of the graph to another vertex $j$ (the value 
$D_{i,j}$) for any $i$ and $j$, we can try to reconstruct the weighted graph
(which can represents a graph in the human body, a hydraulic web...).

Some references on 
the problem of the realization of metric spaces by trees or, 
more generally, by graphs can be found
for instance in \cite{Le} or in the recent paper \cite{D-H-S}, just to quote
 two among  many possible papers.

In this short note,
 we consider a basic graph theory problem, which is, 
in some way, linked to the quoted ones.
Consider a weighted graph $G$ with $n$ vertices, numbered by the set 
$\{1,...,n\}$.  We define the multiset 
$${\cal D}_{i,j} (G) = \{w_G(p) | \; p \; simple \; 
path \; between \; i\, and \;j \}$$
In \S3 we
 establish a criterion to say when, given a multisubset of $ \mathbb{R}$,
there exists a weighted complete graph $G$ 
such that the multisubset is ${\cal D}_{i,j} (G)$ for some $i,j$ vertices of 
$G$ (see Theorem \ref{unsolodij}). 

Besides, in \S4, we establish a criterion to say when, given for
 any $i, j \in \{1,...,n\}$ a multisubset of $ \mathbb{R}$,
${\cal D}_{i,j} $, there is a weighted complete 
graph $G$ such that 
${\cal D}_{i,j} (G)= {\cal D}_{i,j} $ for any $i,j$ (see
 Theorem \ref{complete}).

\section{Notation and remarks}

\begin{notation}  By simple path in a graph, we will mean an unoriented path 
with distinct vertices. If the graph is simple, 
we will denote a path  by the sequence of their vertices.
Obviously $(v_1,.., v_k)$ and $ (v_k,..., v_1) $ are the same path.
A simple path between $i$ and $ j$ ($i$ and $j$ 
vertices of the graph) will denote a simple path whose ends are $i$ and $j$.

 Let $[n]=\{1,...,n\}$.
Let $K_n$ denote, as usual, the complete graph with $[n]$ as set of  vertices.
\end{notation}

\begin{remark} The number of the simple paths between two vertices in $K_n$ 
is $$ N_n := 1 + (n-2) + (n-2) (n-3) +....+ (n-2)!$$
\end{remark}

{\em Proof.} Let $i$ and $j$ be two vertices of $K_n$. 
Obviously there is only one simple path between $i$ and $j$ with only one edge
and there are $n-2$ paths between $i$ and $j$ with $2$ edges.

There are $\left( \begin{array}{c} n-2 \\ 2 \end{array} \right) 2!$ 
simple paths
 between $i$ and $j$ with $3$ edges (we have to choose the two vertices 
of the path besides $i$ and $j$ and their order in the path).

More generally 
there are $\left( \begin{array}{c} n-2 \\ k-1 \end{array} \right) (k-1)!$ 
simple paths
 between $i$ and $j$ with $k$ edges (we have to choose the  $k-1$ vertices 
of the path besides $i$ and $j$ and their order in the path).

Then, in all, they are 
$$ 1 + (n-2) + \left( \begin{array}{c} n-2 \\ 2 \end{array} \right) 2!
+ \left( \begin{array}{c} n-2 \\ 3 \end{array} \right) 3!
+ ........+
 \left( \begin{array}{c} n-2 \\ n-2 \end{array} \right) (n-2)!=$$
$$= 1+ (n-2) + \frac{(n-2)!}{(n-4)!} +  \frac{(n-2)!}{(n-5)!} +....+(n-2)!=$$
$$= 1+ (n-2) + (n-2)(n-3) +....+(n-2)!$$
\hfill {\em Q.e.d.}

\begin{remark} \label{NnNn-1}
$N_n -1 = (n-2) [1 + (n-3) +.... + (n-3)!]= (n-1) N_{n-1}$
\end{remark}

\begin{definition}
Given a weighted simple  graph $G$, for any path $p$ in $G$, we define 
$w_G(p)$ (or simply $w(p)$) 
 as the sum of the weights of the edges of $p$.
 We define the multiset
$${\cal D}_{i,j} (G) = \{w_G (p) | \; p \; simple \; 
path \; between \; i\, and \;j \}$$
\end{definition}

\begin{definition}
Given two sets $S,T \subset \mathbb{R}$ 
of the same cardinality, we define a ``reciprocal
order'' for $S$ and $T$  a bijection $f:S \rightarrow T$. We define 
the difference of $S$ and $T$ by $f$ as the multiset 
$$ \{s-f(s) | \; s \in S \}$$ 
\end{definition}

\begin{definition}
Let $ Y$ be a multisubset of $  \mathbb{R}$ whose elements $y_{l,m}$ 
are indexed by the $2$-subsets $\{l,m\}$
 of $ [n]$ (we write  $y_{l,m}$  instead 
of $y_{\{l,m\}}$). Let $h_{i,j}(Y)$ be the multisubset of 
$  \mathbb{R}$ given by the 
elements $$ y_{i,i_1}+ y_{i_1, i_2} + ... + y_{i_{r-1}, i_r} +y_{i_r, j}$$ 
for $r \in \mathbb{N}$, $i_1,..., i_r \in [n]-\{i,j\}$ distinct. 
\end{definition}

\section{Theorem \ref{unsolodij}}

\begin{remark} \label{rem} 
Let $G$ be a weighted complete graph with $ [n]$ as set of 
vertices and  let $ i,j \in [n]$. Then 

1) for any $ l,m \in [n] - \{ i,j\}$,  
$$ w(l,m) = \frac{1}{2} \left( w(i,l,m,j) + w( i,m,l,j) - w(i,l,j) -w(i,m,j)  \right)$$ 

\begin{center}
\includegraphics[scale=0.4]{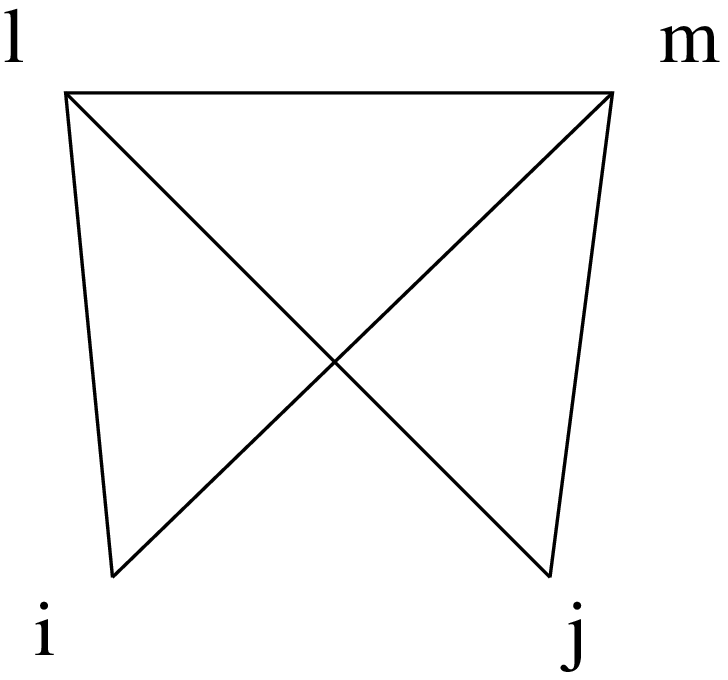}
\end{center}

2) for  any $ l,m \in [n] - \{ i,j\}$,

$ \left\{ \begin{array}{l}
w(i,l)+ w(j,l)= w(i,l,j) \\
w(i,m) + w(j,m) = w(i,m,j) \\
w(i,m) + w(j,l) = w(i,m,l,j) -w(l,m) \\
w(i,l) + w(j,m) = w(i,l,m,j) -w(l,m) 
\end{array} \right.$ 

3) for  any $ i_1,..., i_r \in [n] - \{ i,j\}$, we have that  
$w(i, i_1,...,i_r, j)$ is equal to
$$  \frac{1}{2} \left(   w(i, i_1, i_r, j ) - w(i,
 i_r, i_1, j )+
\sum_{s=1,..., r-1} (w(i, i_s, i_{s+1}, j) + w(i , i_{s+1} ,i_s, j))
\right) -\sum_{s=2,..., r-1} w(i, i_s , j)  $$ 

4) for  any $ l,o,m \in [n] - \{ i,j\}$, 
$$w(i,m,l,j) + w(i,o,m,j) +w(i,l,o,j)=
 w(i,l,m,j) + w(i,m,o,j) + w(i,o,l,j)$$

\end{remark}

{\em Proof.}
The only statement that needs a calculation is 3:
$$ w(i, i_1,...,i_r, j) =  w(i, i_1) + w(i_1, i_2) +... + w(i_{r-1}, i_r) + 
w(i_r, j) =$$
$$= w(i ,i_1, i_r, j) -w(i_1, i_r) +
 w(i_1 ,i_2) +... + w(i_{r-1}, i_r)  = $$ 
$$=\frac{1}{2} \left(   w(i, i_1 ,i_r, j ) -w(i, i_r, i_1, j )+
\sum_{s=1,..., r-1} (w(i, i_s, i_{s+1} ,j) + w(i , i_{s+1}, i_s, j))
\right) -\sum_{s=2,..., r-1} w(i, i_s , j)  $$ 
where the last equality holds by part 1.
\hfill {\em Q.e.d.}

\begin{theorem} \label{unsolodij}
Let $Y$ be  a multisubset of  $ \mathbb{R}$ of cardinality $N_n$.
There exists a weighted 
complete graph $G$ with  $[n]$ as set of vertices 
such that $Y$ is equal to $ {\cal D}_{i,j}(G)$ for some $i,j \in [n]$
 if and only if
we can index the elements $y$ of $Y$ by the finite sequences of elements
in $[n]$ from $i$ to $j$ without repetitions in such a way that:

a) for any $r >2$ and  any $ i_1,..., i_r \in [n] - \{ i,j\}$
$$ y_{i,i_1,...,i_r, j} = \frac{1}{2} \left(   y_{i,i_1,i_r, j } 
- y_{i, i_r, i_1, j }+
\sum_{s=1,..., r-1} (y_{i, i_s, i_{s+1}, j} + y_{i , i_{s+1}, i_s ,j})
\right) -\sum_{s=2,..., r-1} y_{i, i_s,  j}  $$ 

b) for  any $ l,o,m \in [n] - \{ i,j\}$
$$y_{i,m,l,j} + y_{i,o,m,j} +y_{i,l,o,j}=
 y_{i,l,m,j} + y_{i,m,o,j} +y_{i,o,l,j}$$
\end{theorem}

{\em Proof.}
$\Rightarrow$ Follows from Remark \ref{rem}.

$\Leftarrow$ 
Let $G$ be a weighted complete graph with $[n]$ as set of vertices 
and whose weights of the edges are defined in the following way:

$ \bullet $ 
for any $ l,m \in [n] - \{ i,j\}$ we define
 $$ w(l,m) = \frac{1}{2} \left( y_{i,l,m,j} + y_{i,m,l,j} - y_{i,l,j} 
-y_{i,m,j}  \right)$$ 

$ \bullet $ we define $w(i,l)$ and $ w(j,l) $  
for  $ l $ varying in $ [n] - \{ i,j\}$ as  solutions of the linear system $(\ast)$: 

\medskip

$ \left\{ \begin{array}{l}
w(i,l)+ w(j,l)= y_{i,l,j} \\
w(i,m) + w(j,m) = y_{i,m,j} \\
w(i,m) + w(j,l) = y_{i,m,l,j} -w(l,m) \\
w(i,l) + w(j,m) = y_{i,l,m,j} -w(l,m) 
\end{array} \right.$  \hspace{2cm} $ l,m $ varying in $ [n] - \{ i,j\}$ 

\medskip

Observe that the linear system is solvable (but the solutions are not unique), 
in fact:

- the linear system given only by the four above equations with
 $l$ and $m$ fixed
has coefficient matrix 
\begin{center}
\includegraphics[scale=0.6]{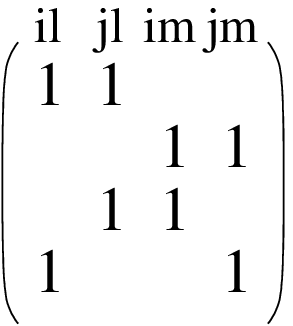}
\end{center}
where the empty  entries are zero and we wrote $il ,\; jl ,\; im ,\;  
jm$, instead of the variables
 $w(i,l)$, $w(j,l)$, $ w(i,m)$,  $ w(j,m)$, over the matrix;
the rank is $3$; precisely the sum of the first two rows is equal to the 
sum of the last two rows; so it is solvable if and only if 
$$y_{i,l,j} + y_{i,m,j} = y_{i,m,l,j} -w(l,m) + y_{i,l,m,j} -w(l,m) $$
 which is true by the definition of $ w(l,m)$

 - the linear system $(\ast)$ has coefficient matrix
\begin{center}
\includegraphics[scale=0.6]{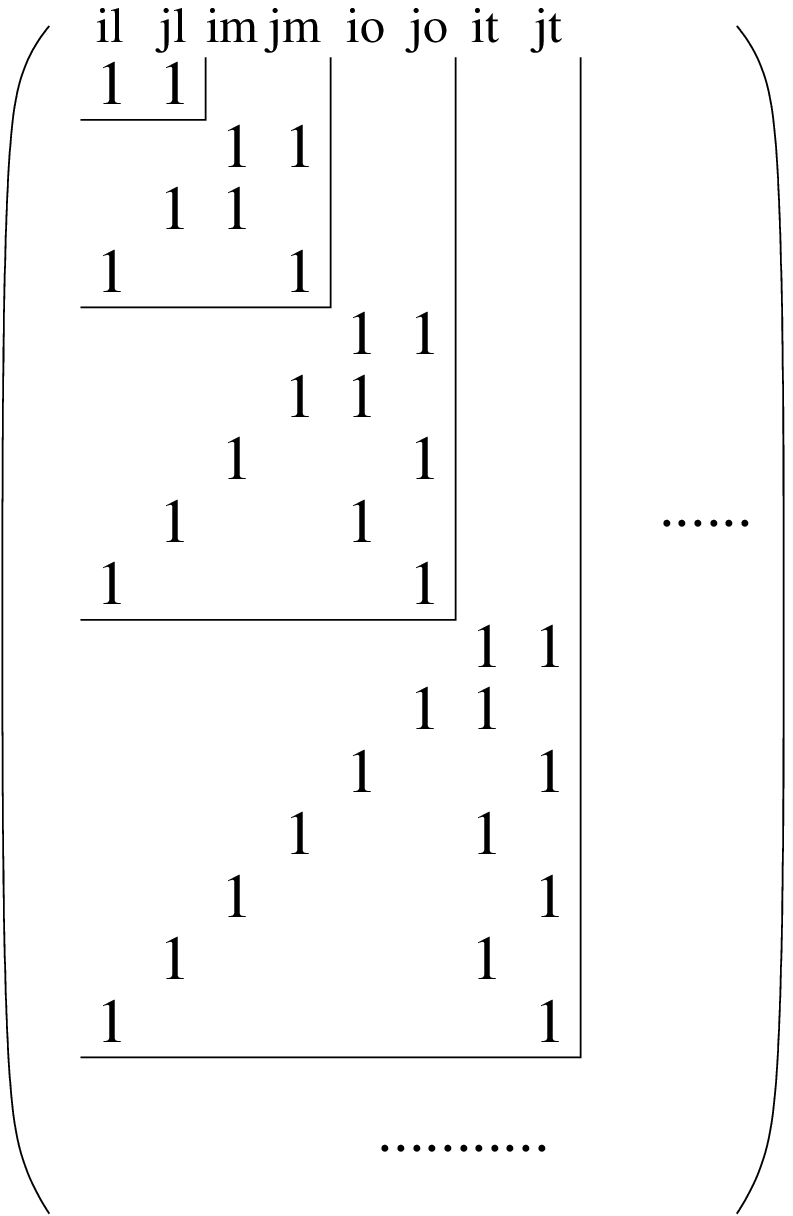}
\end{center}
whose empty entries are zero;
%{\small $$ \left( \begin{array}{ccccccccc}
%il & jl & im & jm &  io & jo & it & jt & ... \\
%\hline
%1  & 1  & 0  & 0  & 0   & 0  & 0  & 0  & ... \\
%0 & 0  & 1 & 1 & 0 & 0 & 0 & 0  & ...\\
%0 & 1  & 1 & 0 & 0 & 0 & 0 & 0  & ...\\
%1 & 0  & 0 & 1 & 0 & 0 & 0 & 0  & ... \\
%0 & 0  & 0 & 0 & 1 & 1 & 0 & 0  & ... \\
%0 & 0  & 0 & 1 & 1 & 0 & 0 & 0  & ... \\
%0 & 0  & 1 & 0 & 0 & 1 & 0 & 0  & ... \\
%0 & 1  & 0 & 0 & 1 & 0 & 0 & 0  & ... \\
%1 & 0  & 0 & 0 & 0 & 1 & 0 & 0  & ... \\
%0 & 0  & 0 & 0 & 0 & 0 & 1 & 1  & ... \\
%0 & 0  & 0 & 0 & 0 & 1 & 1 & 0  & ... \\
%0 & 0  & 0 & 0 & 1 & 0 & 0 & 1  & ... \\
%0 & 0  & 0 & 1 & 0 & 0 & 1 & 0  & ... \\
%0 & 0  & 1 & 0 & 0 & 0 & 0 & 1  & ... \\
%0 & 1  & 0 & 0 & 0 & 0 & 1 & 0  & ... \\
%1 & 0  & 0 & 0 & 0 & 0 & 0 & 1  & ... \\
%  &    &   & .  &  . &  . &   &    &
%\end{array} \right)$$}
in fact the linear system, written explicitely,  is 

\medskip
$ \left\{ \begin{array}{l}
w(i,l)+ w(j,l)= y_{i,l,j} \\
w(i,m) + w(j,m) = y_{i,m,j} \\
w(i,m) + w(j,l) = y_{i,m,l,j} -w(l,m) \\
w(i,l) + w(j,m) = y_{i,l,m,j} -w(l,m) \\ 
w(i,o)+ w(j,o)= y_{i,o,j} \\
w(i,o) + w(j,m) = y_{i,o,m,j}- w(o,m) \\
w(i,m) + w(j,o) = y_{i,m,o,j} -w(o,m) \\
w(i,o) + w(j,l) = y_{i,o,l,j} -w(o,l) \\
w(j,o) + w(i,l) = y_{i,l,o,j} -w(o,l) \\
w(i,t) + w(j,t) = y_{i,t,j} \\
w(i,t) + w(j,o) = y_{i,t,o,j} -w(t,o) \\
w(i,o) + w(j,t) = y_{i,o,t,j} -w(t,o) \\
w(i,t) + w(j,m) = y_{i,t,m,j} -w(t,m) \\
w(j,t) + w(i,m) = y_{i,m,t,j} -w(t,m) \\
w(i,t) + w(j,l) = y_{i,t,l,j} -w(t,l) \\
w(j,t) + w(i,l) = y_{i,l,t,j} -w(t,l) \\
................
\end{array} \right.$ 
\medskip

The coefficient matrix has rank equal to the number of its columns minus $2$;
in fact observe that the rows from the $7^{th}$ to the $9^{th}$, the rows 
from the $12^{th}$ to $ 16^{th}$ and so on are linear combinations of the 
previous ones, for instance the $ 3^{rd}$  plus the $5^{th}$ is equal 
to the sum of the   $ 7^{nd}$ plus the $8^{th}$.

Obviously the linear system is solvable if and only if 
the same relations hold for the constant terms and this is true 
if and only if, for  any $ l,o,m \in [n] - \{ i,j\}$,
$$y_{i,m,l,j}  -w(m,l)+ y_{i,o,j} =  
y_{i,m,o,j} -w(m,o) + y_{i,o,l,j} -w(o,l)$$  
 which is equivalent to assumption b.

\medskip

Now we show that 
$w(i, i_1,...., i_r, j) = y_{i, i_1, .... ,i_r, j}$ for any
 $ i_1,..., i_r \in [n]- \{i,j\}$.

First suppose  that $r=1$. We have that $w(i,i_1, j)= w(i, i_1 ) +w(i_1, j) $  
which is equal to $y_{i, i_1 ,j}$ by the definition of the $w(i,l)$ and 
$w(j,l)$ for $l \in [n] -\{i,j\}$.  
Analogously if $r=2$.

If $r>2$, by assumption a, we have that 
$$ y_{i,i_1,...,i_r, j} = \frac{1}{2} \left(   y_{i,i_1,i_r, j } - y_{i, i_r, i_1, j }+
\sum_{s=1,..., r-1} (y_{i, i_s, i_{s+1}, j} + y_{i , i_{s+1}, i_s ,j})
\right) -\sum_{s=2,..., r-1} y_{i, i_s,  j} $$  which (by cases $r=1,2$)
 is equal to
$$\frac{1}{2} \left(   w(i ,i_1, i_r, j ) - w(i, i_r ,i_1, j )+
\sum_{s=1,..., r-1} (w(i, i_s, i_{s+1}, j) + w(i , i_{s+1}, i_s ,j))
\right) -\sum_{s=2,..., r-1} w(i, i_s,  j)   $$
which is equal to $ w(i, i_1,...,i_r, j)$ by Remark \ref{rem}, part 3.
\hfill {\em Q.e.d.}

\section{Theorem \ref{complete}}

Roughly speaking, the following theorem says that, 
given  multisubsets of $ \mathbb{R}$,  $ {\cal D}_{i,j} $,  
of cardinality $N_n$, for   $i,j \in [n]$,
there exists a weighted complete graph  $G$  such that 
${\cal D}_{i,j} (G) = {\cal D}_{i,j}$  for all $i,j $ if and only if 
we can order the  $ {\cal D}_{i,j} $ reciprocally in such a way that the 
difference of each pair of them can be divided in
into $n-2$  multisubsets, one 
 of cardinality $ N_{n-1} +1$, the others of cardinality $ N_{n-1} $, such
 that all the elements of each  of these subsets have all the same 
absolute value and one of the  $ {\cal D}_{i,j} $
is in the image of $h$.

\begin{theorem} \label{complete}
For any $i,j \in [n]$,  let $ {\cal D}_{i,j} $ be a  multisubset of 
$ \mathbb{R}$ of cardinality $N_n$. There exists a weighted complete
 graph  $G$ with vertices $[n]$ such that 
$${\cal D}_{i,j} (G) = {\cal D}_{i,j}\;\; \forall i,j \in [n]$$ 
if and only if, for any $i,j,k \in [n]$, there exists a reciprocal 
order $f^{i,k}_{j,k}$ for ${\cal D}_{i,k}$ and ${\cal D}_{j,k}$ 
and for any $i,k \in [n]$ an element $y_{i,k} \in {\cal D}_{i,k}$   such that

A)  $f_{j,k}^{i,k} (y_{i,k}) =y_{j,k}$

B) the difference of  ${\cal D}_{i,k}$ and ${\cal D}_{j,k}$  
by   $f_{j,k}^{i,k}$ can be  divided into $n-2$  multisubsets which we call 
$${\cal L}_0 (i,k| j,k), \;\;\; {\cal L}_r (i,k | j,k)$$ 
for $r \in [n] -\{i,j,k\}$,
the first of cardinality $ N_{n-1} +1$, the others of cardinality $ N_{n-1} 
$, such that 
%all the elements of each  of these subsets have  the same absolute value
 ${\cal L}_0 (i,k| j,k)$ contains one element equal to $y_{i,k} -y_{j,k}$ 
and the other elements are equal to its opposite and 
${\cal L}_r (i,k | j,k)$, for $r \in [n] -\{i,j,k\}$,
contains $N_{n-2}$ elements  equal to $y_{r,i} -y_{r,j}$ 
and  the others are equal to its opposite.

C) if  $Y:=\{y_{l,m}\}_{l,m \in [n]}$, 
  there exist $u,v \in [n]$ such that ${\cal D}_{u,v} = h_{u,v} ( Y)$.
\end{theorem}

{\em Proof.}
$\Rightarrow$ Let $ y_{i,j}= w_G(i,j)$. 
We can divide the paths between $i$ and $k$ in two kinds: the ones
passing through $j$, which we can write as $ (i ,\gamma, j ,\eta, k)$ for 
some $\gamma$ and $ \eta$ disjoint subsets  of $[n] -\{i,j,k\}$,
and the others, which we can write as  $(i ,\delta, k) $ with $\delta$ 
subset  of $[n] -\{i,j,k\}$.

We establish a bijection $f=  f^{i,k}_{j,k}$ 
between the paths between $i$ and $k$ and the paths
between $j$ and $k$ which will define a reciprocal order of  
$ {\cal D}_{i,k}$ and $ {\cal D}_{j,k}$:
$$  (i ,\gamma ,j ,\eta, k ) \stackrel{f}{\longmapsto} 
(j ,\gamma^{-1}, i, \eta ,k)$$ 
$$  (i, \delta, k) \stackrel{f}{\longmapsto} (j ,\delta, k)$$
for  $\gamma$ and $ \eta$ disjoint subsets  of $[n] -\{i,j,k\}$,
$\delta$ subset  of $[n] -\{i,j,k\}$ 
(if $ \gamma = (\gamma_1,..., \gamma_t)$, we denote
 $ \gamma^{-1} = (\gamma_t,..., \gamma_1)$).

\begin{center}
\includegraphics[scale=0.3]{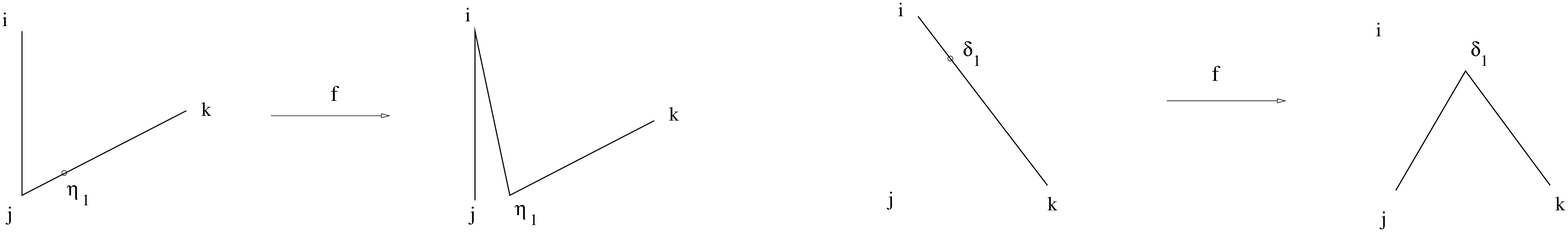}
\end{center}

The paths between $i$ and $k$ can be divided into $n-2$ subsets:

$\bullet $ a subset ${\cal P}_0  $ 
of cardinality $N_{n-1}+1$ whose elements are:

- the path $(i,k)$ 

- the paths of the kind $ (i, \gamma, j ,k)$ with $\gamma \subset  [n]
 -\{i,j,k\}$

(that is the paths with no vertices between $k$ and $i$ or $k$ and $j$ if 
they pass through $j$) 

$\bullet$ $n-3$ subsets ${\cal P}_{r} $, for $ r \in [n] -\{i,j,k\}$,
each of them  defined as follows:  fix 
 $r \in [n] -\{i,j,k\}$ and consider 

- the paths of kind $ (i, \eta, k)$ 
with $ \eta \subset  [n] -\{i,j,k\}$,  $\eta \neq \emptyset$ and
 $\eta_1=r$

- the paths of kind $ (i, \gamma, j ,\eta, k)$  with 
 $\gamma, \eta \subset  [n] -\{i,j,k\}$, $\eta \neq \emptyset$ and $\eta_1=r$

(that is the paths with some vertices between $k$ and $i$ or $k$ and $j$ if 
they pass through $j$). 

Obviously all the subsets  ${\cal P}_{r} $, for 
$ r \in [n] -\{i,j,k\}$,
have the same cardinality, so the cardinality of each of them is
$$\frac{ N_n -1 -N_{n-1}}{n-3}
= \frac{(n-2) N_{n -1} -N_{n-1}}{n-3}= N_{n-1}$$
where the first equality holds by Remark \ref{NnNn-1}.

Observe that for any $ p \in  {\cal P}_0  $ 
we have $$|w(p)- w(f(p))|= |w(i,k)-w(j,k)|= |y_{i,k}-y_{j,k}| $$ 
 in particular 

$ w(i,k)- w(f(i,k))= w(i,k) -w(j,k)= y_{i,k}-y_{j,k}$

$ w(i, \gamma, j, k )- w(f(i, \gamma, j ,k))= 
w(i, \gamma ,j, k )- w(j ,\gamma^{-1} ,i, k))
 = w(j,k) -w(i,k)= -y_{i,k}+y_{j,k}$.

Besides for any  $ p \in  {\cal P}_r $,  , for 
$ r \in [n] -\{i,j,k\}$,
we have $$|w(p)- w(f(p))|= |w(i,r)-w(j,r)| = | y_{i,r}-y_{j,r}|$$ 
 in particular 

$w( i, \eta, k) -w(f( i ,\eta, k))= w( i, \eta, k) -w(j ,\eta, k)=
 w(i, \eta_1) -w(j ,\eta_1)=  w(i ,r) -w(j, r) = y_{i,r}-y_{j,r}$ 

$w(i, \gamma, j ,\eta, k)-w(f(  i ,\gamma, j ,\eta, k))=
w(i, \gamma ,j ,\eta, k)-w(j ,\gamma^{-1}, i, \eta, k))=
 -w(i ,\eta_1) +w(j ,\eta_1)= - w(i ,r) +w(j ,r) = -y_{i,r}+y_{j,r}$. 

So the sets ${\cal P}_0$ and ${\cal P}_r$ give the subsets 
${\cal L}_0$ and ${\cal L}_r$ in 
the difference of  ${\cal D}_{i,k}$ and ${\cal D}_{j,k}$  by  $f^{i,k}_{j,k}$ 
and B  holds (A and C are obvious).

$ \Leftarrow$  Let $G$ be the weighted complete graph with $[n]$ as set 
of vertices and whose weights are defined by
$$ w_G(i,k) = y_{i,k} $$ 
For the graph $G$ we can define a set of  reciprocal orders 
${}^G \! f^{i,k}_{j,k}$ for 
the ${\cal D}_{i,k}(G)$ as in the proof of the other 
implication. By such reciprocal orders 
 the difference of  ${\cal D}_{i,k}(G)$ and ${\cal D}_{j,k}(G)$  
 can be  divided into $n-2$  multisubsets,
${\cal L}^G_0 (i,k| j,k)$ and $  {\cal L}^G_r (i,k | j,k)$
for $r \in [n] -\{i,j,k\}$,
the first of cardinality $ N_{n-1} +1$, the others of cardinality $ N_{n-1} 
$, such that   ${\cal L}^G_0 (i,k| j,k)$ 
contains one element equal to $w_G(i,k) -w_G(j,k)$
 and the other elements are equal to its opposite and 
${\cal L}^G_r (i,k | j,k)$, for $r \in [n] -\{i,j,k\}$,
contains $N_{n-2}$ elements   equal to $w_G(r,i) -w_G(r,j)$
and  the others are equal to its opposite.

We have to show that  $$ {\cal D}_{i,k} (G)= {\cal D}_{i,k}$$
for any $i, k \in [n]$.

First we want to prove that, for any $i,j,k$, the difference of
$ {\cal D}_{i,k}(G)$ and $ D_{j,k}(G) $ by ${}^G \! f^{i,k}_{j,k}$
 is equal to the 
difference of $ {\cal D}_{i,k}$ and $  D_{j,k}$ by $ f^{i,k}_{j,k}$.
Obviously ${\cal L}_0^G (i,k| j,k) ={\cal L}_0 (i,k| j,k)$, because
 ${\cal L}_0^G (i,k| j,k) $ will be composed by 
$  w_G(i,k) -w_G(j,k)  $ 
   and $N_{n-1}$ opposites of it and 
 ${\cal L}_0 (i,k| j,k) $ will be composed by 
$ y_{i,k} -y_{j,k}  $    and $N_{n-1}$ opposites of it by assumption; 
so from our definition of the weights of $G$ we can conclude.
Also ${\cal L}_r^G (i,k|j,k) ={\cal L}_r (i,k|j,k)$ for $r \in [n] -\{i,j,k\}$ 
because ${\cal L}_r^G (i,k|j,k)$ is given by $N_{n-2}$ numbers equal to
$ w_G(i,r)- w_G(j,r) =y_{i,r}- y_{j,r}$ 
and $N_{n-1} -N_{n-2}$ equal to its opposite
and  the same ${\cal L}_r (i,k|j,k)$.

\medskip

 We have that $ {\cal D}_{u,v}(G)=h_{u,v} (\{w_G (l,m)\}_{l,m})= h (\{y_{l,m}
\}_{l,m})=   {\cal D}_{u,v} $, where the last equality holds
 by assumption C. From 
this and from the fact that the difference of
$ {\cal D}_{k,u}(G)$ and $ D_{u,v}(G) $ by $ {}^G f^{k,u}_{u,v}$ is equal to the 
difference of $ {\cal D}_{k,u}$ and $  D_{u,v}$ by $f^{k,u}_{u,v}$,
 we get that 
$ {\cal D}_{k,u}(G)= {\cal D}_{k,u} $ for any $k \in [n] -\{u,v\}$.

 From 
this and from the fact that the difference of
$ {\cal D}_{i,k}(G)$ and $ D_{k,u}(G) $ by $ {}^G \!
 f^{i,k}_{k,u}$ is equal to the 
difference of $ {\cal D}_{i,k}$ and $  D_{k,u}$ by $f^{i,k}_{k,u}$,
 we get that 
$ {\cal D}_{i,k}(G)= {\cal D}_{i,k} $ for any $i,k$. 
\hfill {\em Q.e.d.}

\end{document}